# Parabolic resonances and instabilities.

Vered Rom-Kedar
The Department of Applied Mathematics and Computer Science
The Weizmann Institute of Science,
P.O.B. 26, Rehovot 76100, Israel.
vered@wisdom.weizmann.ac.il

October 3, 1996

### Abstract

A *parabolic resonance* is formed when an integrable two-degrees-of-freedom Hamiltonian system possessing a circle of parabolic fixed points is perturbed. It is proved that its occurrence is generic for one parameter families (co-dimension one phenomenon) of near-integrable, two d.o.f. Hamiltonian systems. Numerical experiments indicate that the motion near a parabolic resonance exhibits a new type of chaotic behavior which includes instabilities in some directions and long trapping times in others. Moreover, in a degenerate case, near a *flat parabolic resonance*, large scale instabilities appear. A model arising from an atmospherical study is shown to exhibit flat parabolic resonance. This supplies a simple mechanism for the transport of particles with *small* (i.e. atmospherically relevant) initial velocities from the vicinity of the equator to high latitudes. A modification of the model which allows the development of atmospherical jets unfolds the degeneracy, yet traces of the flat instabilities are clearly observed.

# Contents

# 1 Introduction

The motion of $n$ particles under conservative forces (e.g. the planets in the solar system, coupled frictionless pendulums, point vortices etc.) is described by an $n$-degrees-of-freedom Hamiltonian system. This motion may be ordered or chaotic, depending on the Hamiltonian $H(x, v)$ (the energy), the dimension of the system and the initial position and velocities of the particles. Consider a Hamiltonian flow for which all solutions are ordered. A fundamental question is what will be the solution's fate if the system is changed just by a small amount (e.g. the effect of a meteor on the solar system). In particular, consider solutions which are bounded and are "roughly stable" in the unperturbed system (there exists a finite size neighborhood in phase space which stays "near" these solutions forever). Can solutions in this neighborhood become unbounded or considerably different if the system is changed by little? In this paper we discover a new low dimensional mechanism ($n = 2$) creating such a dramatic effect on solutions. We prove it is a mechanism which is physically common, and demonstrate its consequences on a simple model for the motion of weather balloons in the atmosphere.

    The appearance of homoclinic loops in integrable Hamiltonian systems provides a mechanism by which by small perturbation integrability is destroyed and new complex dynamics arises. This phenomenon has been prin-

cipally investigated in the context of one-and-a-half d.o.f. Hamiltonian system. Yet, even with regards to this simplest setting for which homoclinic chaos exists, important open questions linger, see for example papers in this issue by Klafter, Meiss, Melnikov and Zaslavsky and references therein.

What is the structure of near-integrable two d.o.f. Hamiltonian system? Generically, most phase space is foliated by two-dimensional tori which divide the energy surfaces. Substantially different dynamics appear near singular points or curves on the energy surface - either fixed points or invariant circles (or lines in the non-compact case) of the flow [15]. When isolated unstable fixed points exist inherently three dimensional phenomena arise, and these are not yet fully understood. See Lerman and Umanskii [15]for the classification of integrable systems and Koltsova and Lerman [16]and Grotta-Ragazzo [14]for the analysis of chaotic system with a homoclinic loop to a saddle-center.

Consider near-integrable Hamiltonians of the form:

$$H(x,v,D,\theta;\epsilon) = H_0(x,v,D)+\epsilon H_1(x,v,D,\theta;\epsilon), \qquad (x,v,D) \in R^3, \theta \in T, \tag{1.0.1}$$

Such a form arises from a general Hamiltonian $H(x,v,y,w;\epsilon)$ of a near integrable flow by using one of the constants of motion of the unperturbed flow ($\epsilon = 0$), $D(x,v,y,w)$ as a variable, and defining $\theta$ as its conjugate variable (e.g. $D$ is the angular momentum). Consider a domain $U$ for which the transformation $(x,v,y,w) \to (x,v,D,\theta)$ is non-singular[1]. In particular, in $U$, a fixed point in the $(x,v)$ plane, $q_f = (x_f, v_f, D_f)$:

$$\nabla_{(x,v)} H_0(x_f, v_f, D_f) = 0 \tag{1.0.2}$$

corresponds to an invariant circle in the four-dimensional phase space (and not to an isolated fixed point). This circle may be composed of a family of periodic orbits (when $\dot\theta(q_f) \neq 0$) or of a circle of fixed points (when $\dot\theta(x_f, v_f, D_f) = 0$). An important lesson from Lerman and Umanskii [15]work is that phenomena which are co-dimension one in near-integrable one-and-a-half d.o.f. Hamiltonian systems may become generic in integrable two d.o.f. Hamiltonian systems; In the $(x,v)$ plane, $q_f$ may generically be hyperbolic or elliptic, and generically, there exists $D_f = D_p$ values for which it is *parabolic* [15]. Similarily, generically, there exists $D_f = D_R$ values for which $\dot{\theta}(x_f, v_f, D_f) = 0$. These simple observations are the basis for the proof that parabolic circles of fixed points (i.e. $D_f = D_R = D_p$) are generic in one-parameter family of integrable two d.o.f. Hamiltonian systems, see section 3.

---

[1]i.e. where the surfaces spanned by the Poisson action of the integrable system [15] are at least one-dimensional

The behavior near the invariant circles under small perturbations depends on their stability (i.e. the stability of $q_f$ in the $(x, v)$ plane), on the rotation rate on them (i.e. $\dot{\theta}(q_f)$), and the rotation rate on orbits which are homoclinic to them (if such exist). If $q_f$ is elliptic, stability reigns, though one may obtain resonant or nonresonant behavior depending on the frequency ratios. In particular, if $q_f$ is an elliptic circle of fixed points, the resonance creates localized structures in $\theta$. When $q_f$ is hyperbolic, and the energy surfaces are compact, under small perturbation several types of homoclinic tangles may appear. It is conjectured that only one type, the reducible homoclinic tangle, is topological (orbital) conjugate to the one appearing in the one-and-a-half d.o.f. Hamiltonian system case.

The *reducible homoclinic tangle* appears when, in the integrable limit, there exists a homoclinic loop to a periodic orbit, and the motion along the homoclinic loop and its fat neighborhood preserves the angular direction of the periodic orbit (i.e. $|\dot{\theta}| > a > 0$ along the periodic orbits and the homoclinic loop). In this case a global Poincaré map transverse to the periodic orbit is well defined in a neighborhood of the homoclinic loop and the flow is conjugate to a two-dimensional symplectic map with a homoclinic tangle. For separable two d.o.f. Hamiltonian systems homoclinic loops to periodic orbits of bounded periods always satisfy this requirement as in Holmes and Marsden [10]. In non-separable systems, it is satisfied in regions in phase space for which one of the angle variables varies monotonically in $t$ along the unperturbed separatrices as in Holmes and Marsden[11]. Even in this simplest setting there are some differences between two d.o.f. Hamiltonian systems and one-and-a-half d.o.f. Hamiltonian systems when one considers behavior of physical ensembles of initial conditions [1, 2] or relations between, for example residence time and escape rates in the Poincaré section. These differences seems to be physically important yet mathematically - technical.

A second type of homoclinic tangle (*lumpy homoclinic tangle*) arises when small perturbations are applied to systems for which non-monotonicity in the angle occurs along the unperturbed homoclinic solutions. Then, non-uniformity in the angle variable of the chaotic zone is created [1]. Beyond the differences in terms of the observed non-uniformity between the first and second types, it is conjectured that this case is dynamically distinct from the reducible homoclinic tangle, at least for small yet finite perturbation parameter.

A third type of homoclinic tangles, the *homoclinic hyperbolic resonance*, occurs when a system containing a normally hyperbolic circle of fixed points with a family of heteroclinic orbits connecting points on the circle is perturbed, see Kovacic and Wiggins [8], Haller [9]and Kaper and Kovacic [7]. Using ideas as in Lerman and Umanskii [15], it is proved in section 3 that the occurrence of a circle of fixed points which is normally hyperbolic whithin the energy surface is generic in the class of integrable two d.o.f. Hamiltonian systems. Thus, if the energy surfaces are compact, the scenario of homoclinic hyperbolic resonance appears generically in near-integrable two d.o.f. Hamiltonian systems.

In this paper a new type of chaotic behavior is introduced - the *homoclinic-parabolic resonance*. Parabolic invariant circles appear generically in integrable non-separable two d.o.f. Hamiltonian system [15]. When there exists periodic motion on the invariant circles the effect of small perturbation on initial conditions in its neighborhood seems to be insignificant (in fact the separatrices splitting is exponentially small in the distance from the bifurcation point [12]). In particular, initial conditions starting close to the invariant circle stay close to it - the invariant circle is "roughly stable". This changes dramatically when the invariant circle is a circle of fixed points in the four dimensional phase space, see figure 2.0.1 - figure 2.0.5. First, it is proved that this event is not so exotic - it is a co-dimension one phenomena for integrable two d.o.f. Hamiltonian system. Then, it is demonstrated numerically that small perturbations lead to dramatically different dynamics than observed in the other types of homoclinic tangles or the other types of resonances. It is a combination of the hyperbolic homoclinic chaos, the localization of the elliptic resonances and instabilities formed by sliding along the newly formed elliptic circles as seen in slow passage through bifurcations [13]. These effects are even more pronounced when additional degeneracy appears - the flat parabolic resonance case. Here, in the unperturbed system the parabolic invariant circle of fixed points is a part of a family of invariant circles of fixed points all of which belonging to the same energy level set $H_0 = h$. Namely, on this energy surface the non-degeneracy conditions for the iso-energetic KAM theory fails to exist and instabilities reign. This may be an example for a

two-d.o.f. equivalent of the stochastic-webs [18] which appear in one-and-a-half d.o.f. Hamiltonian system(global instabilities will appear if non-compact energy surfaces are considered).

This study has been initiated while investigating a simple model for particles (weather balloons or floats) transport in the atmosphere or ocean. With Y. Dvorkin and N. Paldor [1], we have analyzed the simplest possible two-dimensional motion on a rotating sphere, due only to Coriolis force (i.e. in the absence of any body force, see Paldor and Killworth [6]), perturbed by the inclusion of a zonally travelling pressure wave which perturbs the geopotential surface from its mean spherical shape, see Paldor and Boss [5]. In [1] we developed tools for delineating the phase space of non-separable two d.o.f. Hamiltonian system to regions in which the different types of homoclinic structures appear. Here, it is shown that this model exhibits a flat parabolic resonance when the perturbation is of the form of a standing pressure wave. Adding a latitude dependent pressure gradient, as was suggested in [1], removes this degeneracy. If the pressure gradient is comparable to reasonable atmospherical velocities (the atmospherical relevant case) then strong instabilities induced by the nearby flat parabolic resonance are observed.

The parabolic resonance is a low co-dimension phenomenon hence it is expected to appear in many other applications. For example, the motion of a particle in a central field [4, 11] with effective potential with at least one maximum is an integrable two d.o.f. Hamiltonian system, which, when perturbed by small amplitude travelling wave [2] (in the angle direction) may exhibit a parabolic resonance for some isolated wave speeds (this is a variation on the example of Holmes and Marsden [11] who considered this problem in the hyperbolic regime and with standing wave perturbation). Another example of similar form is the whirling pendulum perturbed by a travelling wave perturbation. Notice that these examples (and the atmospherical model) are formally two-and-a-half d.o.f. problems (since time is explicitly introduced in the travelling wave perturbation), but, as is demonstrated in section 2 and discussed in section 3, this class of systems can be easily reduced to autonomous two d.o.f. Hamiltonian system. Finally, the example which motivated the works on the homoclinic hyperbolic resonance (see Kovacic and Wiggins [8]), namely that of a two-mode truncation of the damped and driven Sine-Gordon equation (see McLaughlin paper in this issue and references therein) also exhibits parabolic resonance for some parameter values.

This paper is organized as follows: In section 2 the simple quasi-inertial model for Lagrangian motion in the atmosphere is recalled and studied numerically near its parabolic resonances. In section 3 it is established that the occurrence of a parabolic resonance is generic for near-integrable two d.o.f. Hamiltonian system depending on one parameter. Then several physically typical forms of two d.o.f. Hamiltonian system (such as separable systems) are considered, and their compliance with the genericity assumptions is examined. Section 4 is devoted for discussion.

## 2  Model for particles' motion in the Atmosphere

The motion of a particle on a rotating sphere subject to a conservative travelling wave perturbation serves as a simple model for the motion of particles in the atmosphere [6, 5, 1]. It was recently suggested that the inclusion of latitude dependent pressure gradient may be incorporated into this model to simulate the appearance of jets in the atmosphere [1]. Using polar co-ordinates, the non-dimensional Lagrangian momentum equations for the eastward and northward velocity components $(u, v)$ and the rate of change of the longitude and latitude coordinates $(\lambda, \phi)$ in the presence of a zonally travelling pressure wave and a latitude dependent pressure gradient are given by:

$$\begin{aligned}
\frac{d\lambda}{dt} &= \frac{u}{\cos\phi} \\
\frac{du}{dt} &= v\sin\phi(1 + \frac{u}{\cos\phi}) - k\epsilon\frac{A(\phi)}{\cos\phi}\cos(k\lambda - \sigma t). \\
\frac{d\phi}{dt} &= v \quad (2.0.3)\\
\frac{dv}{dt} &= -u\sin\phi(1 + \frac{u}{\cos\phi}) - B'(\phi) - \epsilon A'(\phi)\sin(k\lambda - \sigma t) \\
&\phi \in [-\frac{\pi}{2}, \frac{\pi}{2}], \ \lambda \in [0, 2\pi], \ u, v \in \mathrm{R}^2.
\end{aligned}$$

$B(\phi)$ and $A(\phi)$ represent the latitude dependent amplitudes of the pressure gradient and the pressure wave respectively. For simplicity, both are assumed to be even in $\phi$. For $\epsilon = 0$, the system is autonomous, $\lambda$ is uncoupled and

(2.0.3) has two integrals of motion, corresponding to (twice) the angular momentum $D$ and the energy $E$:

$$D = \cos\phi(\cos\phi + 2u) \qquad (2.0.4)$$
$$E = \frac{1}{2}(u^2 + v^2) + B(\phi). \qquad (2.0.5)$$

In [1] the structure of (2.0.3) is studied for $B(\phi) = 0$ for general values of the parameters. By constructing the colored energy surfaces $H = h$ and the corresponding bifurcation diagrams of the energy-momentum map (indicating singular curves, their stability and the monotonicity properties of $\theta(t)$ in the $(H,D)$ plane), the phase space is delineated to regions in which inherently three-dimensional dynamics appears and to regions in which the system is reducible to one-and-a-half d.o.f. Hamiltonian system [1]. Critical parameter values for which the dynamics changes qualitatively are found by this analysis. Specifically, the case of standing waves ($\sigma = 0$) has been identified as a degenerate case. It is shown below that for $B(\phi) = 0, \sigma = 0$, a flat parabolic resonance occurs, and that by appropriate choice of $B(\phi)$ the resonance may be made non-degenerate. Moreover, the behavior near elliptic and hyperbolic resonances is compared with that appearing near a flat and generic parabolic resonances.

Consider the case $k \neq 0$, and define:

$$\theta = \frac{1}{2}(\lambda - ct), \qquad c = \frac{\sigma}{k}. \qquad (2.0.6)$$

Then (2.0.3) becomes:

$$\frac{d\theta}{dt} = \frac{1}{4}\frac{D}{\cos^2\phi} - \frac{1}{2}(c + \frac{1}{2})$$
$$\frac{dD}{dt} = -2k\epsilon A(\phi)\cos(2k\theta)$$
$$\frac{d\phi}{dt} = v$$
$$\frac{dv}{dt} = \frac{1}{8}\sin(2\phi)(1 - \frac{D^2}{\cos^4\phi}) - B'(\phi) - \epsilon A'(\phi)\sin(2k\theta)$$

which is an autonomous system with Hamiltonian

$$H(\phi, v, \theta, D; k, c, \epsilon) = \frac{v^2}{2} + \frac{1}{8}(\frac{D}{\cos\phi} - \cos\phi)^2 - \frac{c}{2}D + B(\phi) + \epsilon A(\phi)\sin(2k\theta)$$

$$k \neq 0, \ -\frac{1}{2} \leq c, \epsilon > 0, \text{sign}(k) = \text{sign}(c). \tag{2.0.7}$$

Since $A(\phi)$ and $B(\phi)$ are assumed to be even[2], the cylinder $\phi = v = 0$ is invariant for all $\epsilon$. On it the equations for $(D, \theta)$ correspond to the pendulum equations with $2k$ islands of width $\sqrt{8\epsilon A(0)}$ centered at $D = D_r(c) = 1 + 2c$, and rotational orbits filling the rest of the cylinder [1]. The behavior near this cylinder in the four-dimensional phase space depends on its stability properties in its normal ($(\phi, v)$) direction.

Consider the unperturbed flow in the $(\phi, v)$ plane. The origin is hyperbolic for $|D| < \sqrt{1 - 4B''(0)}$, having two symmetrical homoclinic orbits extending to $\pm(\phi_h, 0)$ ($\phi_h \approx \arccos(D)$ for $B(\phi) \ll 1$). It is parabolic at $D = D_p = \sqrt{1 - 4B''(0)}$ and elliptic for $|D| > D_p$. Parabolic resonance arises when $D_p \approx D_r(c)$, namely for $c$ values near:

$$c_p = \frac{1}{2}(D_p - 1) = \frac{1}{2}(\sqrt{1 - 4B''(0)} - 1). \tag{2.0.8}$$

For $|D| < D_p$ there are two elliptic fixed points $(\pm\phi_{ell}, 0)$ at

$$D_g^{\pm}(\phi_{ell}) = \cos^2(\phi_{ell})\sqrt{1 - \frac{8B'(\phi_{ell})}{\sin(2\phi_{ell})}}, \tag{2.0.9}$$

with the natural frequency in the $\phi, v$ plane of nearly $\sin \phi_{ell}$ (for small $B'(\phi), B''(\phi)$), and with the frequency in the longitude direction of

$$\frac{d\theta}{dt}\Big|_{(\phi_{ell},0,\theta,D_g(\phi_{ell}))} = \frac{1}{4}\sqrt{1 - \frac{8B'(\phi_{ell})}{\sin(2\phi_{ell})}} - \frac{1}{2}\left(c + \frac{1}{2}\right) \tag{2.0.10}$$

Notice that for $B'(\phi) \ll 1$

$$\frac{d\theta}{dt}\Big|_{(\phi_{ell},0,\theta,D_g(\phi_{ell}))} = \frac{B'(\phi_{ell})}{\sin(2\phi_{ell})} - \frac{1}{2}c. \tag{2.0.11}$$

It follows that when $B(\phi_{ell}) \equiv 0$ (more generally if $B'(\phi) = \alpha \sin(2\phi)$) and $c = c_p = 0$ (resp. $c = c_p \approx 2\alpha$) a degenerate situation occurs by which $\frac{d\theta}{dt}|_{(\phi_{ell},0,\theta,D_g(\phi_{ell}))} = -\frac{1}{2}c_p = 0$. Namely, for each $\phi_{ell} \in (-\pi/2, \pi/2)$, the circle $(\phi_{ell}, 0, \theta, D_g^+(\phi_{ell}))$ consists of fixed points in the four dimensional phase space. Moreover, for $\epsilon = 0$ all these circles of fixed points belong to the same energy surface $H = 0$. In section 3 it is established that this *flat parabolic resonance* is not generic.

---

[2]Thus this assumption of symmetry is convenient, yet it seems to be inessential for what follows.

Numerical simulations of (2.0.7) are performed using DSTOOL [17] with the Bulirsch-Stoer integrator. The Hamiltonian value is monitored so that at least 5 significant digits are unaltered in each computation. In all figures $k = 3$, $A(\phi) = \cos^3(\phi)$, $B(\phi) = \frac{1}{3}\beta \cos^3(\phi)$, and the other parameters, including $\beta = -B''(0)$, are varied as indicated in the captions.

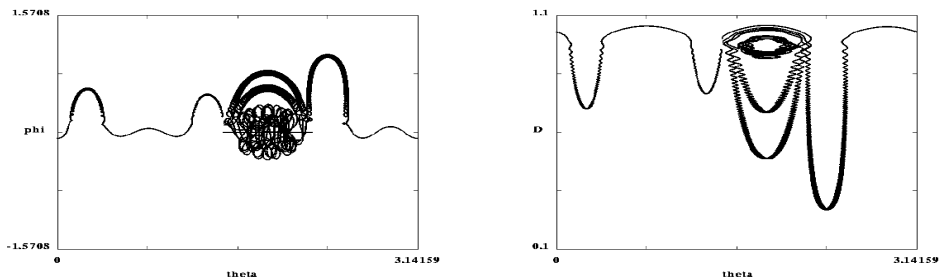

Figure 2.0.1: Flat parabolic resonance.
One trajectory shown in the a) $(\phi, \theta)$ and b) $(D, \theta)$ planes.
$\epsilon = 0.00025, c = 0, \|(\phi_0, v_0, u_0)\| = O(10^{-5}), B(\phi) \equiv 0$.

figure 2.0.1 demonstrates that strong instabilities appear near a flat parabolic resonance. In this figure one trajectory starting close to the flat parabolic resonance at $D_0 = 1$ is shown. The trajectory which is positioned initially close to the equator with very small initial velocities (10 cm/second $= u_0 = 10^{-4}$) travels to latitudes as high as $60^o$ (similar behavior is observed for initial westward velocities of up to 1m/second). Moreover, the motion is restricted for a long period of time to one cell of the resonance band (zonally localized) and then jumps to a different cell. In the $(\theta, D)$ plane it is seen that the motion for $D > D_p$ follows closely the resonant motion of the invariant cylinder $\phi = v = 0$ whereas the instabilities are created for $D < D_p$ when the motion descends along resonant elliptic points up to very small $D$ values ($D = 0.1$ is seen in the figure) creating the instabilities that are seen in the $(\phi, \theta)$ plane. As $-u_0$ is increased (above $-u_0 \approx 0.01$) the descendent along these elliptic points disappears, and as it is increased further (above $-u_0 \approx 0.05$) localized structures disappear. Therefore, a mechanism for the transport of particles launched near the equator with $||(u_0, v_0)|| \ll 1$ to high altitudes has been found, a phenomena that has been observed in atmospherical experiments [3]. Curiously, it appears that the instabilities develop only for sufficiently *small* initial velocities. Typical initial conditions seem to be trapped in one zone for a long periods of time (e.g. figure 2.0.2).

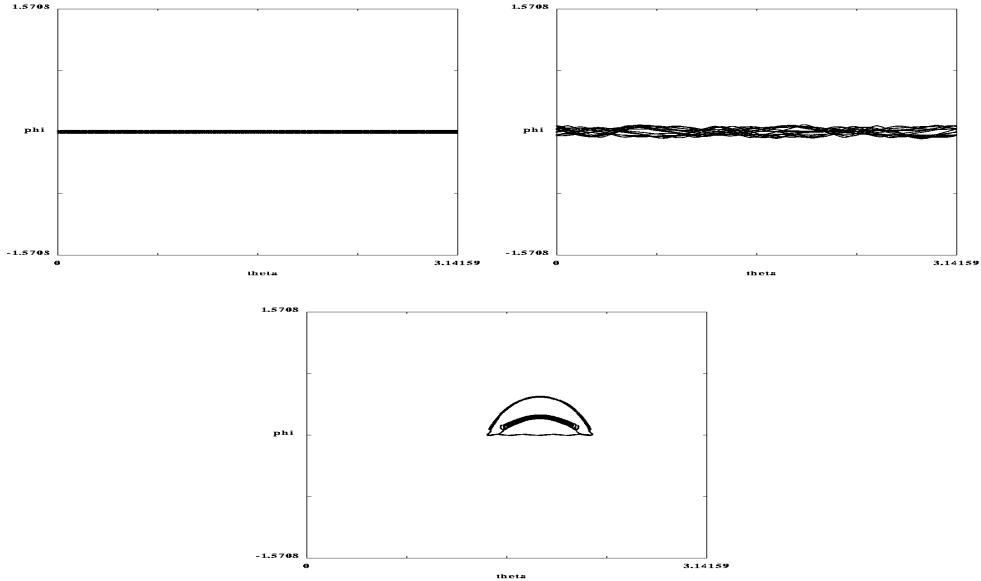

Figure 2.0.2: Near a flat resonance.
The same initial condition is evolved with different wave speeds:
$\epsilon = 10^{-5}, B(\phi) \equiv 0, (\phi_0, v_0, \theta_0, u_0) = (10^{-5}, 10^{-5}, 2.2, \frac{1}{2}10^{-11})$.
a) $c = 0.1$, b) $c = 0.01$, c) $c = 0.0001$.

In figure 2.0.2 the behavior near the parabolic resonance is examined as the wave velocity $c$ is varied: the same initial condition (with $D_0 = 1$) is integrated three times with $c = 0.1, 0.01$ and $c = 0.0001$ respectively. It is seen that for $c$ values sufficiently far from 0 the equator vicinity is stable whereas when $c = 0.0001$ the instability and the localization of the nearby flat resonance is clearly seen.

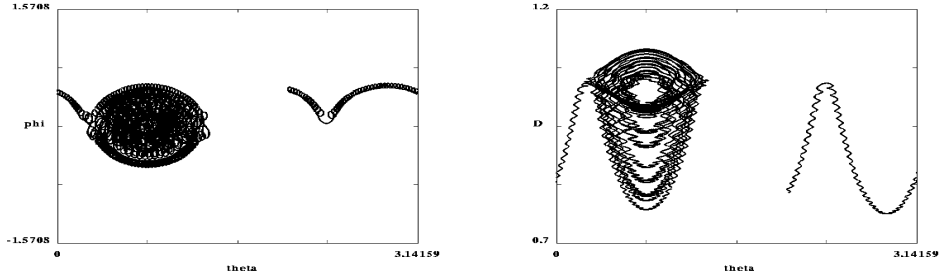

Figure 2.0.3: Nearly flat parabolic resonance.
One trajectory shown in the a) $(\phi, \theta)$ and b) $(D, \theta)$ planes.
$\epsilon = 0.00025, c = 0.029, B''(0) = -0.03, (\phi_0, v_0) = O(10^{-5}), u_0 = 0.00042$.

If $\frac{B'(\phi_{ell})}{sin(2\phi_{ell})} \neq const$ then the flat degeneracy is removed. The effect of removing the flat degeneracy is examined by varying $\beta = B''(0)$. In figure 2.0.3 the nearly flat behavior is examined, with $\beta = 0.03$ ($c_p = 0.0291, D_p = 1.058$), showing that the flat instabilities prevail.

Figures 2.0.4-2.0.5 show the non-degenerate case with $\beta = 0.3$ near $D_p = \sqrt{2.2} \approx 1.48$ and $c = c_p \approx 0.24$. In figure 2.0.4 a trajectory starting very close to the invariant cylinder exhibits strong zonal localization and a mixed behavior of motion near the elliptic points, near the separatrices and near the separatrices of the resonance zone. figure 2.0.5 shows the behavior near the resonance when both $c$ is slightly shifted from $c_p$ and the initial condition of the trajectory is not extremely close to the invariant cylinder as in the previous figure. The chaotic trajectory is presented together with two other trajectories which are on the invariant cylinder ($\phi_0 = v_0 = 0$) to illustrate the shift of the localized structure from its base.

We would like to contrast the behavior observed in these figures, with that observed near elliptic resonance and near hyperbolic resonance; In the elliptic-resonance case near the resonance zone the $D$ variable perform excursions of $O(\sqrt{\epsilon})$. However, in the $(\phi, v)$ plane these excursions are nearly unnoticeable - the motion continues to exist mostly on invariant tori. Some of these tori correspond to oscillatory motion in $\theta$, and small changes in $(\phi, v)$ usually do not alter the oscillatory/rotational behavior in $\theta$, unless the initial conditions are in the extreme (exponentially small?) close vicinity of the separatrices created by the resonance. Thus, in the elliptic-resonance case one observes either regular zonally localized motion or regular motion travelling around the globe and no transitions between the two type of motions is detected by simple numerical experiments (for small $\epsilon$).

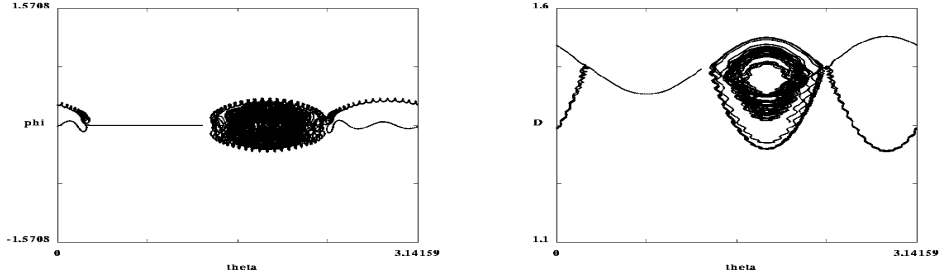

Figure 2.0.4: Parabolic resonance.
One trajectory shown in the a) $(\phi, \theta)$ and b) $(D, \theta)$ planes.
$\epsilon = 0.00025, c = 0.24, B''(0) = -0.3, (\phi_0, v_0, \theta_0, u_0) = (10^{-12}, 10^{-12}, 0, 0.24)$.

The behavior in the vicinity of a hyperbolic-resonance is mostly influenced by the instability in the hyperbolic direction. Thus, the $O(\sqrt{\epsilon})$ oscillations in $D$ are coupled to the unstable chaotic motion in the $(\phi, v)$ plane as is shown in figure 2.0.6 and figure 2.0.7. The extend in latitude to which chaotic trajectories reach is approximately the same as the latitude reached by an unperturbed orbit starting with the same initial velocity $\phi_{\max}^{\epsilon} = \phi_{\max}^{0} + O(\sqrt{\epsilon})$ (see figures). The perturbed motion near the invariant cylinder is chaotic w.r.t. north-south motion and it has a zonal structure scarred by the resonance on the invariant cylinder (see figure 2.0.6 and figure 2.0.7). Localized motion in $\theta$ exists on the invariant cylinder inside the resonance zone. Applying the general theory of Haller [9] for the case $B(\phi) = 0, k \geq 3$, it follows that for $c \in (c_k = -\sin^2(\frac{\pi}{4k}), 0)$ homoclinic orbits to oscillatory periodic orbits exist [1], i.e. that localized homoclinic solutions exist (and may be detected numerically by devicing the correct scheme [8, 9, 7]). Since $c_k \approx -0.066$, this range of $c$ values is also influenced by the flat parabolic resonance, and so the behavior described above for the near-parabolic-resonance regime fits this regime as well. It is unclear yet if this is true in the general case. Taking smaller $c$ values ($c = -0.25$), one finds that trajectories progress along the longitude coordinate (with known "speed" [1, 9]), thus no zonally localized motion is detected outside the invariant cylinder $\phi = v = 0$.

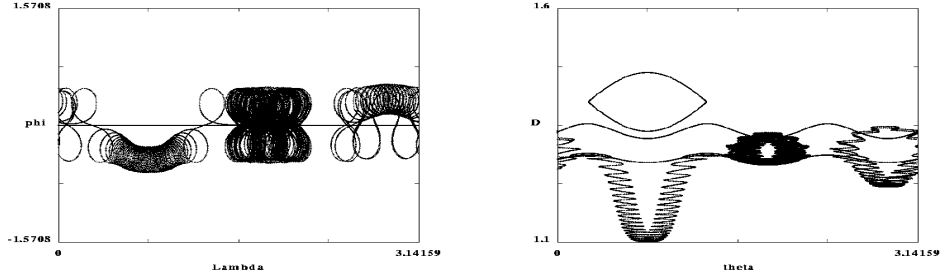

Figure 2.0.5: Near a parabolic resonance.
Three trajectories are shown in the a) $(\phi, \theta)$ and b) $(D, \theta)$ planes.
$\epsilon = 0.00025, c = 0.2(c_p = 0.24), B''(0) = -0.3$,
Two trajectories are on the invariant cylinder $\phi = v = 0$
Third trajectory has $(\phi_0, v_0, \theta_0, u_0) \approx (10^{-5}, 10^{-5}, 0, 0.2)$.

Generic parabolic resonance combines the behavior of the hyperbolic and elliptic regimes - one observes the localized structures due to the elliptic nature of the trajectory for $D > D_p$ and then the hyperbolic nature which induces instabilities eventually breaking this localized structures for $D < D_p$. Moreover, another dominant mechanism for the instability is induced by the elliptic orbits created at $D = D_p$, similar to trajectories passing slowly through a pitchfork bifurcation [13]. Combining these three mechanisms one obtains a motion which has chaotic nature w.r.t. to zonal localization (long trapping periods and then a jump to a different island) and instabilities in the latitude extent of the trajectories. Initial conditions which in the unperturbed case remain in the vicinity of the equator (starting with $D > D_p$) may reach latitudes as high as $25°$ in the non-degenerate case (figure 2.0.4) and much higher latitudes in the flat and nearly flat cases (figure 2.0.1 and figure 2.0.3). Changing the parameters to values which are close to those producing parabolic resonances, show that much of the properties described above are still present - though the zonal localization occurs for lower $D$ values as the distance from the invariant cylinder is increased (figure 2.0.5).

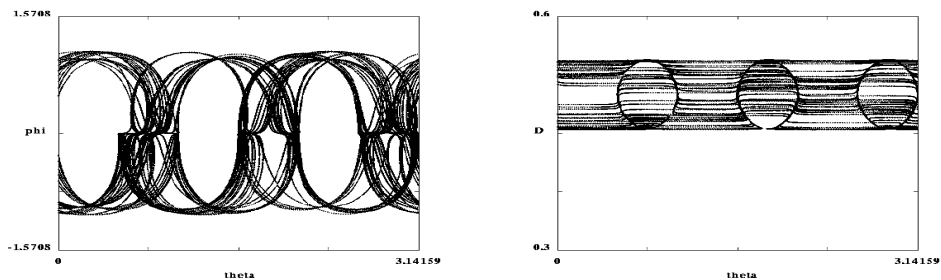

Figure 2.0.6: Hyperbolic resonance, $\theta_0 = 0$.
One trajectory shown in the a) $(\phi, \theta)$ and b) $(D, \theta)$ planes.
$\epsilon = 0.00025, c = -0.25, B(\phi) \equiv 0, (\phi_0, v_0) = O(10^{-5}), u_0 = -0.25$.

## 3  Genericity of parabolic resonances.

Consider a near-integrable two-degrees-of-freedom Hamiltonian flow depending on a parameter $c$:

$$H(x, v, D, \theta; c) = H_0(x, v, D; c) + \epsilon H_1(x, v, D, \theta; \epsilon), \quad (3.0.12)$$
$$(x, v, D) \in R^3, \ \theta \in T, \ c \in R.$$

Such a form arises from a general Hamiltonian $H(x, v, y, w; c, \epsilon)$ of a near integrable flow by using one of the constants of motion of the unperturbed flow, $D(x, p)$ as a variable, and defining $\theta$ as its conjugate variable. Consider a domain $U$ for which the transformation $(x, v, y, w) \to (x, v, D, \theta)$ is non-singular[3]. In particular, in $U$, a fixed point in the $(x, v)$ plane corresponds to either a periodic orbit or to a circle of fixed points in the four-dimensional phase space (and not to an isolated fixed point).

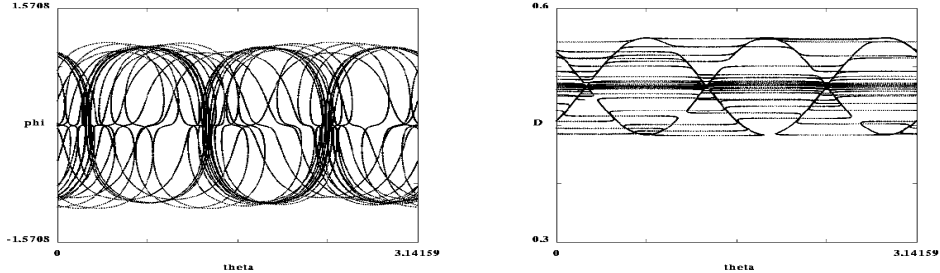

Figure 2.0.7: Hyperbolic resonance, $\theta_0 = 1.8$.
One trajectory shown in the a) $(\phi, \theta)$ and b) $(D, \theta)$ planes.
$\epsilon = 0.00025, c = -0.25, B(\phi) \equiv 0, (\phi_0, v_0) = O(10^{-5}), u_0 = -0.25$.

Generally, one would expect parabolic resonances to be of at least co-dimension 2 in this setting (e.g. one condition for parabolicity and another for vanishing frequency on the invariant circle). Loosely, in view of Lerman and Umanskii [15]work, one parameter may be replaced by the constant of the unperturbed motion, $D$.

It is proved below that transverse intersections of various manifolds in the four dimensional space $(x, v, D, c)$ imply the genericity of the hyperbolic (resp. parabolic) resonances in near-integrable, (resp. one parameter family of near-integrable) two d.o.f. Hamiltonian systems.

Surfaces of fixed points of the two-dimensional flow in the $(x, v)$ plane are determined by:

$$\nabla_{(x,v)} H_0(q_f) = 0, \quad q_f = (x_f, v_f, D_f, c_f), \qquad (3.0.13)$$

hence these are two dimensional surfaces. Denote the union of these surface by $F$ and denote points on $F$ by $q_f$. Generically the fixed points $q_f \in F$ are

---

[3]i.e. where the surfaces spanned by the Poisson action of the integrable system [15] are at least one-dimensional

either elliptic or hyperbolic. Parabolic points $q_p$ on $F$ are determined by the additional condition

$$\det \frac{(\partial^2 H_0(x, v, D; c))}{(\partial x, \partial v)}|_{q_p} = 0 \tag{3.0.14}$$

hence correspond to a one dimensional curve $P \subset F$.

Define the three-dimensional surface $T$ consisting of points $q_T$ at which $\frac{d\theta}{dt} = 0$:

$$\frac{\partial H_0}{\partial D}|_{q_T} = 0 \tag{3.0.15}$$

Generically, $T$ intersects $F$ transversely along a one dimensional curve $R$ of fixed points $q_R \in R$ for which $\frac{d\theta}{dt} = 0$, namely $q_R$ correspond to a circle of fixed point in the full phase space. Generically the fixed point on $R$ are composed of intervals of elliptic fixed points and intervals of hyperbolic fixed points. Since generically,

$$\frac{\partial H_0}{\partial c}|_{q_R} \neq 0 \tag{3.0.16}$$

the curve $R$ intersect the three dimensional surface $c = const$ transversely at a point(or several points), hence, we have established:

**Lemma 1:** *The occurrence of a normally hyperbolic (resp. elliptic) circle of fixed points is generic in integrable two d.o.f. Hamiltonian system.*

**Corollary 2:** *The occurrence of a hyperbolic (elliptic) resonance is generic in near-integrable two d.o.f. Hamiltonian systems.*

Moreover, generically the three dimensional surface $T$ intersects the one dimensional curve of parabolic fixed points $P$ at isolated points $q_{PR}$, which are circles of parabolic fixed points in the four dimensional space. Thus, it has been established that:

**Lemma 3:** *The occurrence of a parabolic circle of fixed points is generic in a one-parameter family of integrable two d.o.f. Hamiltonian systems.*

**Corollary 4:** *The occurrence of a parabolic resonance is generic in a one-parameter family of near-integrable two d.o.f. Hamiltonian systems.*

Finally, notice that the surface of constant energy containing the parabolic circle of fixed points $H_0(q) = H_0(q_{PR})$ defines a three dimensional surface $H_{PR}$, which generically intersects the one dimensional curve of circles of fixed points, $R$, at isolated points, and the two-dimensional surface $F$ along a curve $F_{PR}$. Hence, generically, on its energy surface, $q_{PR}$ is an isolated circle of fixed points belonging to a one dimensional family of invariant circles ($F_{PR}$) on which $\frac{d\theta}{dt} \neq 0$.

The *tangential parabolic resonance* occurs when the intersection of $H_{PR}$ with $R$ is not transverse at $q_{PR}$, namely $F_{PR}$ and $R$ are tangent at $q_{PR}$. Clearly this is a co-dimension two phenomena. Then, there is an (infinitesimal) band of circles of fixed points on $H_{PR}$. A *flat parabolic resonance* appears when the two curves $R$ and $F_{PR}$ coincide. If along $F_{PR}$ $c$ remains constant, then a whole line of invariant circles of fixed points belongs to the same energy surface for the same parameter value - a degenerate situation in which no barriers to transport exist under small perturbations. The atmospherical model described in section 2, exhibits, when $B(\phi) \equiv 0$, this degenerate flat parabolic resonance.

The form (3.0.12) arises naturally in applications, four of which were mentioned in the introduction. In the applications, the Hamiltonian $H_0(x, v, D; c)$ may appear in some special form. The implications on the non-degeneracy conditions mentioned above are examined for three physically common cases:

## 3.1 The separable case.

The separable case arises naturally when two one-d.o.f.-systems are weakly coupled. Then, the transformation $(x, v, y, w) \to (x, v, D, \theta)$ is applied near invariant circle in the $(y, w)$ plane, namely $(D, \theta)$ are the action-angle coordinates of the $(y, w)$ unperturbed flow. Moreover, it is natural in this case to assume that the parameter $c$ appears only in one of the one d.o.f. Hamiltonian systems. Hence we assume

$$H_0(x, v, D; c) = H_0^s(x, v; c_1) + H_0^D(D; c_2), \qquad (3.1.1)$$

or, more explicitly

$$\nabla_{x,v} \frac{\partial H_0}{\partial D} \equiv 0, \qquad (3.1.2)$$

and consider the two cases $c_1 = const$ and $c_2 = const$ separately.

**Fixing** $c_2$, i.e. letting
$$\frac{\partial^2 H_0}{\partial D \partial c} \equiv 0, \tag{3.1.3}$$
does not alter the previous results; Note that the curves corresponding to parabolic invariant circles, $P$, (which generically occur for isolated $c_1$ values), are just straight lines of the form $\{q_p = (x_p, v_p, c_{1p}, D), D \in [D_1, D_2]\}$. Since in this case:
$$T = \{q_T : \frac{\partial H_0^D(D;c_2)}{\partial D}|_{q_T} = 0\} = \{q_T : D_T = D_R\} \tag{3.1.4}$$
it follows that $P$ intersects $T$ transversely. Similarly, the curve $R$ lies in $T$ and generically intersects the surface $c = const$ transversely. It follows that the previous Lemmas and corollaries (1-4) are valid with the words "two d.o.f. Hamiltonian system" replaced by "two d.o.f. Hamiltonian systems satisfying (3.1.2) and (3.1.3)".

**Fixing** $c_1$, corresponds to:
$$\nabla_{x,v} \frac{\partial H_0}{\partial c} \equiv 0 \tag{3.1.5}$$
and produces different behavior; In this case, generically, parabolic invariant circles do not appear. The two dimensional surface of fixed points $F$, is parallel to the $(D, c)$ plane, and on it there exists the curve $R$ on which circles of fixed points live. Thus for this case Lemma 1 and Corollary 2 apply yet Lemma 3 and Corollary 4 are not valid; in fact:

**Lemma 5:** A necessary condition for the generic existence of parabolic circle of fixed points in a one parameter family of integrable two d.o.f. Hamiltonian system is
$$||\nabla_{x,v} \frac{\partial H_0}{\partial c}|| + ||\nabla_{x,v} \frac{\partial H_0}{\partial D}|| \not\equiv 0. \tag{3.1.6}$$

## 3.2 The travelling wave (TW) case.

Consider unperturbed Hamiltonians of the form:
$$H_0(x, v, D; c) = H_0^T(x, v, D) + H_0^D(D; c), \tag{3.2.1}$$

namely,
$$\nabla_{x,v}\frac{\partial H_0}{\partial c} \equiv 0, \nabla_{x,v}\frac{\partial H_0}{\partial D} \neq 0. \qquad (3.2.2)$$

These are called travelling wave Hamiltonians by the following motivation; consider an integrable two d.o.f. Hamiltonian system with canonical coordinates $(x, v, D, \lambda)$, which is perturbed by a travelling wave perturbation in $\lambda$:
$$H(x, v, D, \lambda) = H_0(x, v, D) + \epsilon H_1(x, v, D, k\lambda - \sigma t). \qquad (3.2.3)$$

Formally this Hamiltonian corresponds to a two-and-a-half d.o.f. Hamiltonian system. Using the trivial change of variables (for $k \neq 0$)
$$\theta = \lambda - \frac{\sigma}{k}t = \lambda - ct \qquad (3.2.4)$$

it is reduced to an autonomous two d.o.f. Hamiltonian system of the (TW) form with:
$$\bar{H}(x, v, D, \theta) = H_0(x, v, D) - cD + \epsilon H_1(x, v, D, k\theta). \qquad (3.2.5)$$

Indeed,
$$\dot{\theta} = \dot{\lambda} - c = \frac{\partial \bar{H}}{\partial D} \qquad (3.2.6)$$
$$\dot{D} = -\epsilon k \frac{\partial H_1(x, v, D, \alpha)}{\partial \alpha} = -\frac{\partial \bar{H}}{\partial \theta}. \qquad (3.2.7)$$
$$(3.2.8)$$

Three of the four examples mentioned in the introduction are of this form.

In this case the parabolic invariant circles occur for isolated fixed $D$ values independent of $c$, so the curves $P$ are straight lines parallel to the $c$ axis. The curves $R$ of circles of fixed points are as in the general case, hence Lemma 1-4 are valid for the TW case as well.

## 3.3 The natural mechanical case.

If $H_0$ is constructed from a family of one-degree-of-freedom systems with an effective potential which depends on $D$:
$$H_0(x, v, D; c) = \frac{1}{2}v^2 + V(x, D; c) \qquad (3.3.1)$$

then $q_{PR} = (x_{PR}, 0, D_{PR}, c_{PR})$, and the conditions (3.0.13), (3.0.14), (3.0.15) become:

$$\begin{aligned}
\frac{\partial V}{\partial x}|_{q_{PR}} &= 0 \\
\frac{\partial^2 V}{\partial x^2}|_{q_{PR}} &= 0 \\
\frac{\partial V}{\partial D}|_{q_{PR}} &= 0
\end{aligned} \tag{3.3.2}$$

Generically this form obeys the assumed transversality of the manifolds as in the general case, hence Lemma 1-5 apply. For example, the following Hamiltonian (with $q_{PR} = (0, 0, 1, 0)$) satisfies the genericity assumptions (for $b \neq 0$):

$$H_0(x, v, D; c) = (\frac{1}{2a_3} + b)\frac{1}{2}(D-1)^2 - cD + \frac{1}{2}v^2 + \frac{1}{2}a_1 x^2(1-D) + \frac{1}{3}a_2 x^3 + \frac{1}{4}a_3 x^4 \tag{3.3.3}$$

Taking $a_3 > 0$, $b > -\frac{1}{2a_3}$ corresponds to compact energy surfaces $H_0(x, v, D; c) = h$ for any finite $h$. Adding the perturbation $H_1 = \epsilon(1 - \frac{1}{2}x^2)\cos(k\theta)$, the corresponding vector field is of the form:

$$\begin{aligned}
\dot{x} &= v & (3.3.4) \\
\dot{v} &= a_1(D-1)x + a_2 x^2 + a_3 x^3 + \epsilon x \cos(k\theta) & (3.3.5) \\
\dot{\theta} &= (\frac{1}{2a_3} + b)(D-1) - c - \frac{1}{2}a_1 x^2 & (3.3.6) \\
\dot{D} &= k\epsilon(1 - \frac{1}{2}x^2)sin(k\theta), & (3.3.7)
\end{aligned}$$

attaining hyperbolic resonances for $c > 0$ at:

$$D_r(c) = 1 + \frac{c}{\frac{1}{2a_3} + b} \tag{3.3.8}$$

a parabolic resonance at $D_p = 1$, $c_p = 0$ and a flat parabolic resonance at these values for $b = 0$. Numerical experiments with this model produce similar results to the ones presented for the atmospherical problem described in section 2.

# 4  Conclusions

The behavior near parabolic resonances provides a new type of chaotic mechanism for two d.o.f. Hamiltonian system. It has been established that this phenomenon is of co-dimension 1, hence it is expected to appear in numerous applications. Moreover, the basic idea of Lerman and Umanskii [15], that the addition of another non-separable d.o.f. may be considered from bifurcation theory point of view as the addition of another parameter, applies to higher degrees of freedom as well - e.g. we may think of the wave speed $c$ as another conserved quantity of a larger system showing that for higher dimensional systems the parabolic resonance case may be generic.

Preliminary numerical observations of the behavior near parabolic and flat parabolic resonances are presented in figures 2.0.1 - 2.0.5. The features of the parabolic resonance seems to mix three different dynamical behaviors - localization as in the elliptic resonance case, homoclinic chaos which destroys this localization and slow passage through a bifurcation which leads to strong instabilities.

For the atmospherical model with $(B(\phi) = 0)$, for any given wave speed $-0.5 < c < 0$, a hyperbolic resonance occurs: particles sent from the equator's vicinity with westward velocity $u_0 \approx c(< 0)$ participate in a complex dynamical behavior moving north and south chaotically, reaching latitudes of $\arccos(1 - 2|u_0|) + O(\sqrt{\epsilon})$. The motion of such particles may be non-uniform zonally [1] (i.e. some zones will be visited more frequently than others), yet generally, zonally localized motion is not robust. When the speed of the travelling pressure wave vanishes (i.e. the perturbation is of the form of a standing wave) parabolic circles of fixed points occurs for $\epsilon = 0$. Small perturbations of the form of standing waves or travelling waves with small wave speeds results in parabolic resonance. Moreover, this case is degenerate leading to strong instabilities as shown in figure 2.0.1. This degenerate model thus supplies a mechanism for the transport of particles with small initial velocities near the equator to high latitudes. Adding a latitude dependent pressure gradient $(B(\phi) \neq 0)$, which physically corresponds to including the influence of jets [1] removes the degeneracy. If this gradient is not large (which seems to be a reasonable assumption from the physical point of view as the pressure gradient is proportional to the particles velocity, which should be of order 0.01) the degenerate instabilities still appear (figure 2.0.3).

The results presented are preliminary - analytical and extensive numerical studies are under progress. It appears that near a generic parabolic resonance the instability zone in the latitude coordinate corresponds approximately to the maximal latitude achieved by an unperturbed homoclinic orbit to a hyperbolic circle with $D = D_p - O(\sqrt{\epsilon})$. The behavior near a flat parabolic resonance is not resolved yet even numerically. In particular, the scaling of the radius of instability with respect to the perturbation strength (if it is finite), and with respect to the distance in phase-space and parameter space from the flat parabolic resonance case are yet to be found. General conditions for the occurrence of a flat parabolic resonance are yet to be derived.

Many questions are left unanswered at this point, and in particular the relations of this work to McLaughlin, Klafter, Meiss, Melnikov and Zaslavsky presentations are yet to be explored.

# Acknowledgement

It is a pleasure to thank Y. Gutkin, L. Lerman and D. Turaev for stimulating and inspiring discussions and comments. This research was supported by MINERVA Foundation, Munich/Germany.

# References

[1] V. Rom-Kedar, Y. Dvorkin and N. Paldor [1996] Chaotic Hamiltonian dynamics of particle's horizontal motion in the atmosphere, preprint. – [1995] *Chaotic motion on a rotating sphere*, Levy flights and related topics in Physics Springer verlag lecture notes in Physics, Proc., eds. M.F. Shlesinger, G.M. Zaslavsky and U. Frisch, pg. 72-87.

[2] V. Rom-Kedar [1996] Some characteristics of two-degrees-of-freedom Hamiltonian flows, Marseilles.

[3] Jullian, P., V. Lally, W. Kellogg, V. Suomi and C. Cote [1977] The TWERLE Experiment. Bull. Am. Met. Soc. 58(9), 936-948.

[4] L.D. Landau and E.M. Lifshitz, *Mechanics*, Course of theoretical physics, vol. 1, Pergamon Press, 1960.

[5] N. Paldor and E. Boss [1992] Chaotic Trajectories of Tidally Perturbed Inertial Oscillations, J. Atm. Sc., Vol 49, 23, pp 2306-2318.

[6] N. Paldor and P.D. Killworth [1988] Inertial Trajectories on a Rotating Earth, J. Atm. Sc., Vol. , No. 24, pp 4013-4019.

[7] T. Kaper and G. Kovacic [1994] Multi-bump orbits homoclinic to resonance bands, preprint.

[8] G. Kovacic and S. Wiggins [1992] Orbits Homoclinic to resonances, with an application to chaos in a model of the forced and damped sine-Gordon equation, Physica D, 57 pp 185-225.

[9] G. Haller [1995] Universal homoclinic bifurcations and chaos near double resonances, preprint, G. Haller and S. Wiggins [1995] N-pulse homoclinic orbits in perturbations of resonant Hamiltonian systems, Arch. of rational Mech., 130, 25-101.

[10] P.J. Holmes and J.E. Marsden [1982]. Horseshoes in Perturbations of Hamiltonian Systems with Two Degrees of Freedom, Commun. Math. Phys. 82, 523-544.

[11] P.J. Holmes and J.E. Marsden [1983] Horseshoes and Arnold diffusion for Hamiltonian systems on Lie groups. Indiana University Math. Journal 32, No 2, 273-309.

[12] P. Holmes, J. Marsden and J. Scheurle, Exponentially small Splittings of Separatices with applications to KAM theory and Degenerate Bifurcations, Cont. Math. **81**, 213-244 (1988).

[13] N.R. Lebovitz and A.I. Pesci, Dynamics bifurcation in Hamiltonian systems with one degree of freedom [1995] SIAM J. Appl. Math., Vol 55, 4, 1117-1133. G.J.M Mareé, Slow passage through a pitchfork bifurcation [1996], SIAM J. Appl. Math., to appear.

[14] C. Grotta-Ragazzo [1995] Dynamics and diffusion near double-resonances in Hamiltonian systems with discrete symmetries, preprint.

[15] L.M. Lerman and Ia, L. Umanskii [1992,1993] Classification of four dimensional integrable Hamiltonian systems and Poisson actions of $R^2$ in

extended neighborhood of simple singular points, I,II. Russian Acad. Sci Sb. Math, Vol 77 (94), No. 2, Vol 78, (94) No 2.

[16] O. Yu. Koltsova and L.M. Lerman [1995] Periodic and homoclinic orbits in two-parameter unfolding of a Hamiltonian system with a homoclinic orbit to a saddle-center. Bifurcation and Chaos, Vol 5, no 3.

[17] DSTOOLS, Computer program, Cornell university, center of applied mathematics, Ithaca, NY 14853.

[18] G.M. Zaslavsky, R. Z. Sagdeev, D.A. Usikov and A.A. Chernikov, [1991] Weak Chaos ans Quasiregular Patterns, Cambridge Univ. Press, Cambridge.